\documentclass{notices}
\usepackage{amsfonts,amssymb,amsmath,amscd}

\usepackage{esint}
\usepackage{amsthm,color}
\usepackage{enumerate}
\usepackage{mathtools}
\usepackage{enumitem}
\usepackage{url}
\usepackage{xcolor} 
\usepackage{cancel}
\usepackage{caption,subcaption}
\usepackage{tikz}
\usepackage{tikz-3dplot}

\usepackage{longtable}

\newcommand{\K}{{\mathbb K}}
\newcommand{\R}{{\mathbb R}}
\newcommand{\C}{{\mathbb C}}

\newcommand{\N}{{\mathbb N}}

\newcommand{\cE}{{\mathcal E}}

\newcommand{\cA}{{\mathcal A}}

\newcommand{\cF}{{\mathcal F}}

\def\u{{\mathbf u}}
\def\g{{\mathbf g}}
\def\v{{\mathbf v}}

\def\vphi{{\boldsymbol\phi}}

\def\g{{\mathbf g}}

\def\0{{\mathbf 0}}

\def \mb{\mathbb}

\newcommand{\e}{\varepsilon}

\newcommand{\dist}{\operatorname{dist}}

\newcommand{\Hess}{\operatorname{Hess}}

\theoremstyle{plain}
\newtheorem{thm}{Theorem}

\newcommand{\thistheoremnames}{}
\newtheorem*{genericthms}{\thistheoremnames}
\newenvironment{para*}[1]
  {\renewcommand{\thistheoremnames}{#1}%
   \begin{genericthms}}
  {\end{genericthms}}

\theoremstyle{remark}

\newtheorem*{claim*}{Claim}

\numberwithin{equation}{section}

\title{Constraint Maps and Free Boundaries}

\author{Alessio Figalli
  \affil{
    Alessio Figalli is a professor at Department of Mathematics, ETH Z\"urich. His email address is alessio.figalli@math.ethz.ch.
    }
  \and
  Andr\'e Guerra
  \affil{
    Andr\'e Guerra is a Junior Fellow at the Institute for Theoretical Studies, ETH Z\"urich. His email address is andre.guerra@eth-its.ethz.ch.
    }
  \and
  Sunghan Kim
  \affil{
    Sunghan Kim is a postdoc at Department of Mathematics, Uppsala University. His email address is sunghan.kim@math.uu.se.
    }
\and
Henrik Shahgholian
  \affil{
   Henrik Shahgholian is a professor at Department of Mathematics, KTH Royal Institute of Technology. His email address is henriksh@kth.se.
   }
}

\begin{document}

\maketitle
\setcounter{tocdepth}{1}

\section{A brief history of free boundaries}
Free boundary problems emerge in various fields of partial differential equations. These problems occur when the behavior of variables changes abruptly at certain values. 
Numerous examples abound, such as the solid-liquid interface during material solidification, the boundary between exercise and continuation regions of a financial instrument (like an American option), and the transition from elastic to plastic behavior when stress surpasses a critical threshold. These applications often lead to free boundary problems being referred to as phase transition problems.

Classically, mathematical modeling  of free boundaries emerged in studying minimal surfaces under constraints, and can be traced back to 
J.D. Gergonne \cite{Gergonne}, who in 1816 proposed the following problem:

\begin{quote}
    \textit{Couper un cube en deux parties, de telle manière que la section vienne se terminer
aux diagonales inverses de deux faces opposées, et que l'aire de cette section,
terminée à la surface du cube, soit un minimum. 
}
\end{quote}
In short, this says: ``Cut a cube into two parts so that the section ends at opposite diagonals of two opposite faces, and minimize the area of this section.''
This problem remained untouched for almost half a century, until Schwarz 
in 1872 took on studying it, see \cite{Schwarz}.

\begin{figure}
    \centering
    \includegraphics[width=0.5\linewidth]{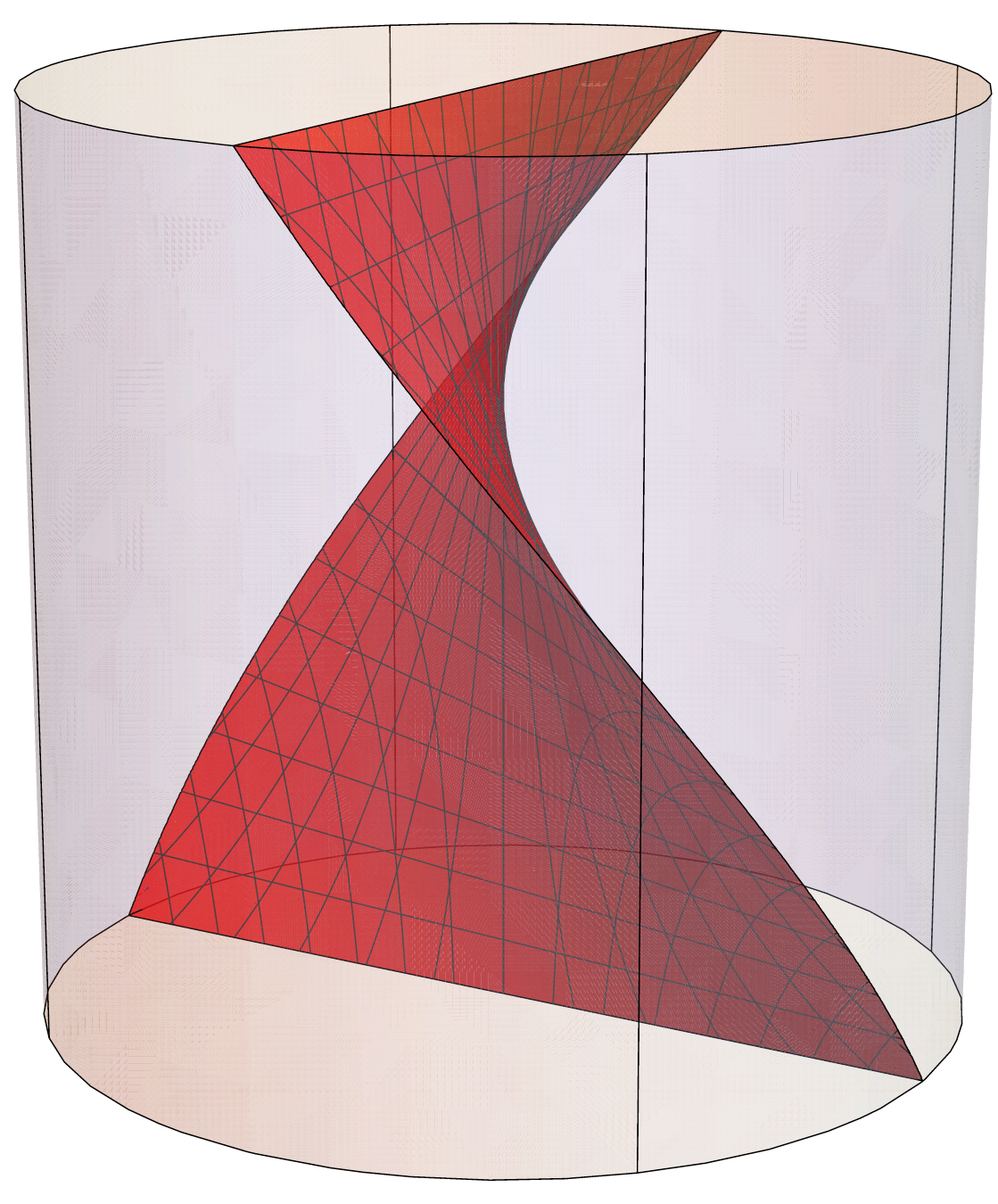}
    \caption {\small The solution to Gergonne's problem in the case where we cut a cylinder, rather than a cube, in half.} \label{fig:Gergonne}
\end{figure}

In fact, in practical applications, the above type of ideas dates back to Roman times, when ancient architects and engineers utilized principles of minimal surfaces to construct stable and efficient structures like arches and domes. These applications focused on optimizing materials and forms to achieve structural stability, indirectly relating to minimizing surfaces under constraints.

In this survey we will focus on a specific but remarkably ubiquitous family of free boundary problems, those of \textbf{obstacle-type}. We will recount the historical development of the topic together with more recent advances, focusing in particular on vector-valued variational problems with constraints. Our discussion is shaped by the contributions of several distinguished mathematicians working at the intersection of free boundary problems, harmonic maps, and geometric partial differential equations.

\section{The scalar obstacle problem}

The obstacle problem is a prominent and perhaps the most widely studied type of free boundary problem. Its theoretical advancements have led to progress in numerous other problems in both applied sciences and theoretical research.

To frame the obstacle problem in a mathematical context, consider a simple 
model of the deformation of an $n$-dimensional membrane in $\R^{n+1}$, keeping in mind the physical description for $n=2$.  We prescribe the deformation of the membrane on its boundary, and we require the membrane to lie above a given obstacle.  For simplicity, we consider the idealized situation in which the membrane is a homogeneous plate of infinite thickness, which does not have resistance, and in which the only acting force is tension\footnote{In general, the membrane may also  be loaded by an external force, which we ignore here.}. We also assume that the deformation of the membrane is prescribed on its boundary.
The solution of the obstacle problem then describes the equilibrium configuration of the membrane under the above constraints.

Let us now give a mathematical formulation of the above problem. 
Denote by
\begin{equation}
    \label{eq:graph}
    x_{n+1}=u(x), \quad x=(x_1, \cdots, x_n)\in\overline{\Omega},
\end{equation}
the deformation representing the equilibrium configuration of the membrane, which we assume to be represented as the graph of a function $u$ defined over the closure of a bounded smooth domain $\Omega\subset \R^n$. Since the deformation is prescribed on the boundary,  we have
\begin{equation}
u(x)=g(x),\ \  x\in\partial\Omega,\label{eq:nd-boundary}
\end{equation}
for some given function $g\colon \partial \Omega\to\R$, which we assume to be smooth.
For an ideal membrane such as the one we described,  the potential energy $\mathcal{P}$ of the deformation
is proportional to the increase in area of the surface of the
membrane, i.e.
$$
\mathcal{P}=  (S_1-S_0),
$$
where $S_0,S_1$ are the surface areas of the membrane before and
after the deformation, respectively; for simplicity, we have set the constant of proportionality to be one.
We thus  have that
$S_0= \hbox{volume}(\Omega)$ and, from basic Calculus, 
\begin{equation*}\label{eq:membrane}
S_1=\int_\Omega \sqrt{1+|D u(x)|^2} \ dx.
\end{equation*}
From the principle of least action, the equilibrium configuration  should minimize the potential energy (i.e., it should minimize the above integral); thus, in the absence of any constraint,  we see that the membrane is deformed to become a \textbf{minimal graph}. If in addition we suppose that the deformation of the membrane is smooth and very small,  then 
$$
S_1 \approx \int_\Omega
\left( 1+\frac{1}{2} |D u(x)|^2 \right) dx.
$$
Therefore, locally,  the potential energy of deformation will essentially be $\tfrac 1 2 \int_\Omega |D u(x)|^2 dx$. Thus, for simplicity, we define
$$
\mathcal{E} (u)= \int_\Omega |D u(x)|^2 dx,
$$
which is the integral studied by Dirichlet, Riemann and many others in classical function theory. Minimizers of $\mathcal{E}$ are \textbf{harmonic functions}, i.e., they satisfy
\begin{equation}
    \label{eq:harmonic}
    \Delta u = 0.
\end{equation}

\begin{figure}
\centering
    \includegraphics[width=5cm]{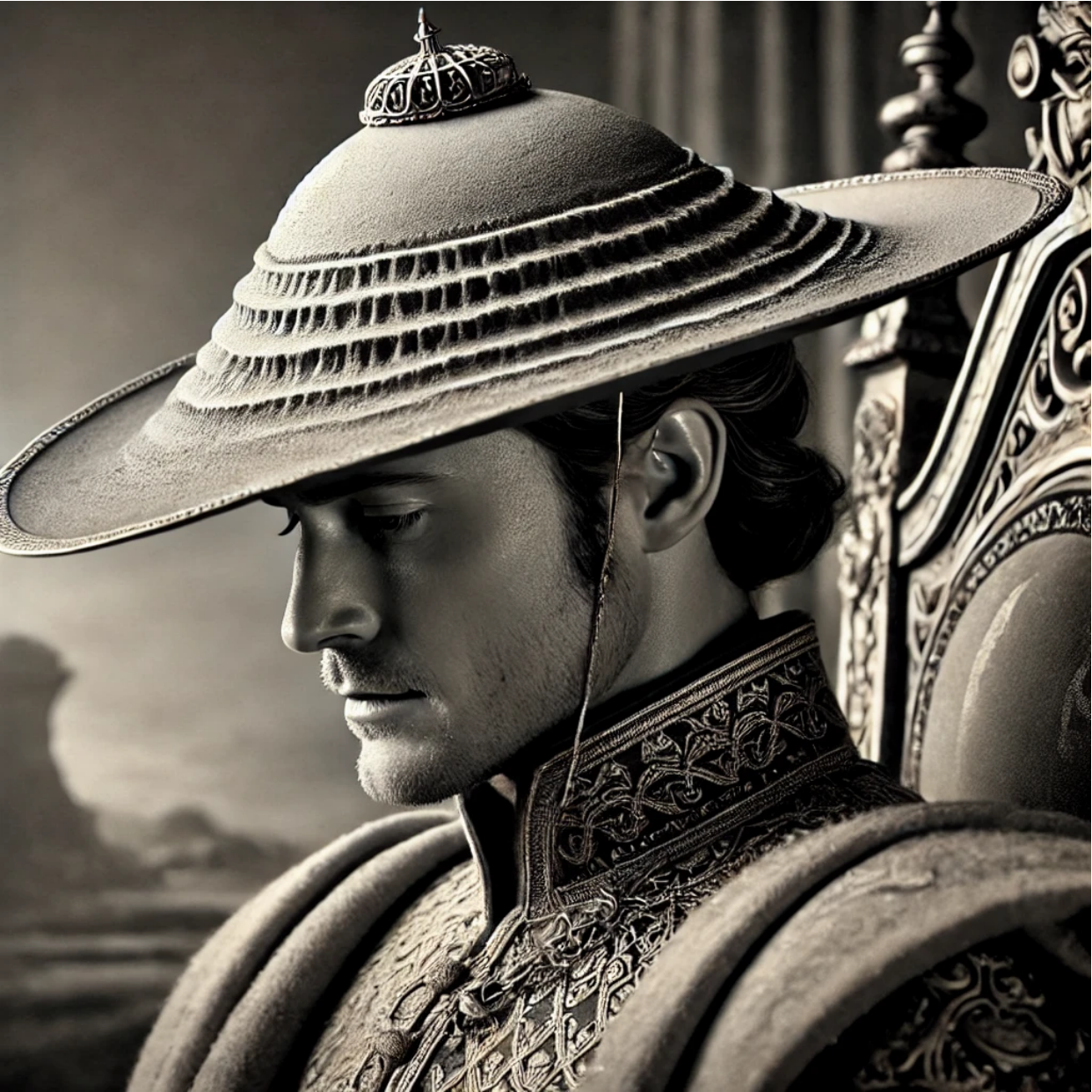}
    \caption {\small AI-generated  King's hat.}
    \label{Fig:king}
\end{figure}

We are yet to incorporate the constraint that the membrane lies above a given obstacle.  An  anecdotal example  goes as follows:
\begin{quote}
{\it The King ordered a hat to be made for himself, fashioned as a minimal garment over a specified ridge.}    
\end{quote}
A hat like this would resemble an example of a \textit{constrained} minimal surface, as illustrated in Figure \ref{Fig:king}. Our task would be to find the configuration among all those in \eqref{eq:nd-boundary} that minimizes the energy $\mathcal{E}$ subject to the obstacle constraint (the surface should lie above the head of the King).

In order to formulate the constraint mathematically, consider an open set $ O\subset  \R^{n+1}$, which we call an
\textbf{obstacle}\index{obstacle!n-dimensional}, and which is of the form
$$
O=\{(x,x_{n+1}): \
x \in\overline{\Omega},\  x_{n+1}<
\psi(x)\},
$$
where $\psi:\overline{\Omega}\to\R$ is a given
smooth function. Our assumption is that the membrane does not penetrate the obstacle, i.e., $(x,u(x))\not \in O$ or, equivalently, due to our graphical assumptions, that the membrane stays above the obstacle:
\begin{equation}
    \label{eq:abovemembrane}
    u(x)\geq \psi(x),\quad x\in\overline\Omega.
\end{equation}
Thus, combining \eqref{eq:nd-boundary} and \eqref{eq:abovemembrane}, we define the class of admissible deformations to be 
$$
\K=\{v\in W^{1,2}_g(\Omega):\  v(x)\geq \psi(x) \text{ for } x\in
\Omega\},
$$
where $ W_g^{1,2}(\Omega)$ is the space of functions in $L^2 (\Omega)$, whose derivatives are also in $L^2(\Omega)$, and which take the boundary value $g$. We remark that $\K$ may be empty, for instance if $g(x_0) < \psi(x_0)$ for some $x_0 \in \partial \Omega$: in this case the problem becomes trivial. 

\begin{figure}
    \centering
    \includegraphics[width=\linewidth]{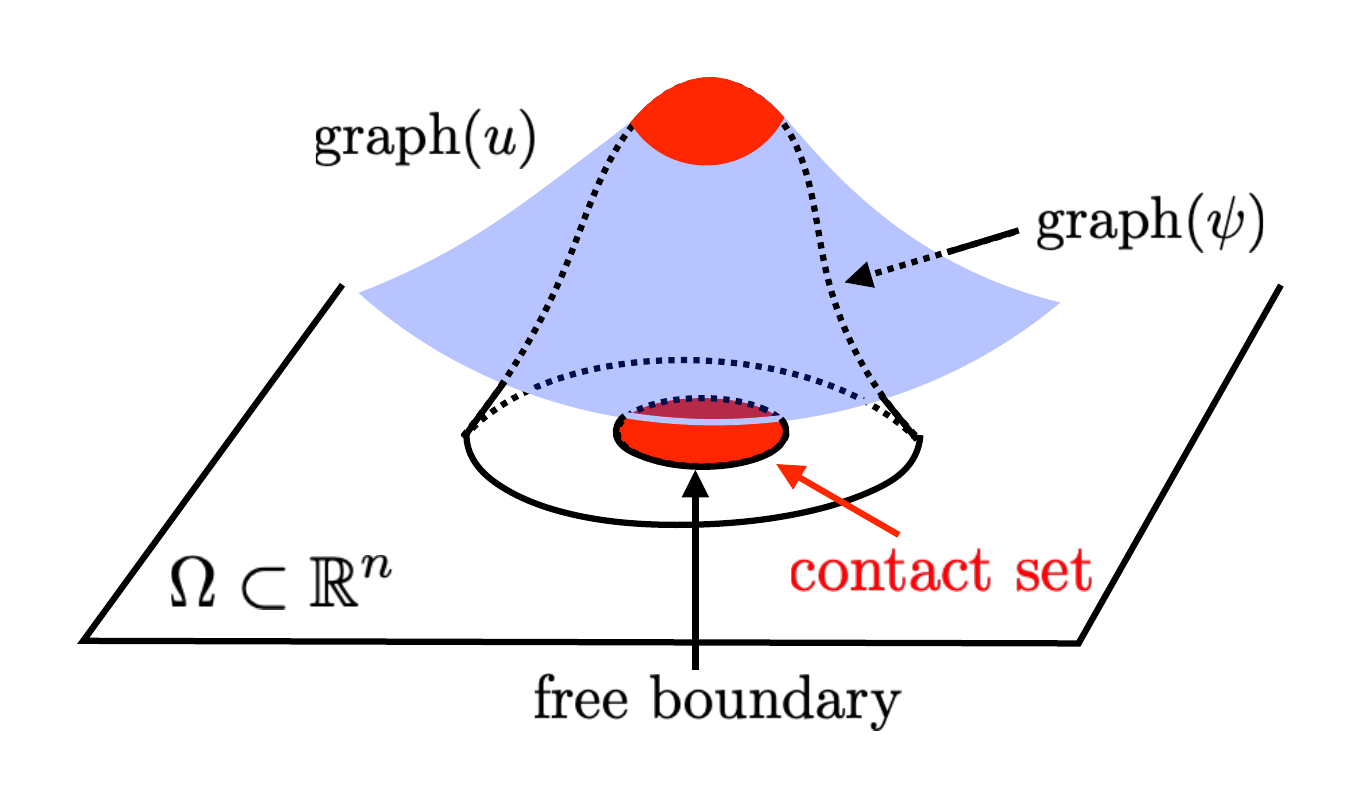}
    \caption {\small A possible solution to the obstacle problem} \label{fig:obstacle}
\end{figure}

In short, the obstacle problem is to find the solution $u\in \K$ to the following problem:
$$\
\mathcal{E}(u)=\min_{v\in \K} \mathcal{E}(v).
$$
One can then show that the solution is characterized through the following partial differential equation
\begin{equation}
    \label{eq:scalarobstacle}
    \Delta u = \Delta\psi\,  \chi_{\{u=\psi\}}.
\end{equation}
This equation asserts that when $u$ does not feel the obstacle (i.e., $u>\psi$) then $u$ is harmonic, as in \eqref{eq:harmonic}; and that when $u$ feels the obstacle (i.e., $u=\psi$) we have $\Delta u = \Delta \psi$. 
Hence \eqref{eq:scalarobstacle} implies that the equation for $\Delta u(x)$ changes depending on which of the following two sets the point $x$ belongs to:
$$\textbf{non-contact set }\{u>\psi\}, \quad \textbf{contact set }\{u=\psi\}.$$
A fundamental feature of the problem is that these sets are not prescribed a priori, and neither is their topological boundary $\partial \{u>\psi\}$, which is therefore called a \textbf{free boundary}. See Figure \ref{fig:obstacle} for an illustration.

Provided that $\K$ is non-empty, standard arguments in the Calculus of Variations yield the existence of a generalized solution to the obstacle problem, and so a basic issue is to understand the regularity and geometric properties of such solutions and of their respective free boundaries. This question is typically considered under the assumption that
\begin{equation}
\label{eq:Deltapsi=1}
\Delta \psi <0.
\end{equation}
In the absence of this condition, the free boundary can behave essentially in an arbitrary way \cite{Caffarelli1998}.
Regarding the regularity of the solution, it has been known since the 1970s that $u \in C^{1,1}(\Omega)$, which is optimal, as \eqref{eq:scalarobstacle} and \eqref{eq:Deltapsi=1}  imply that the Laplacian of $u$ is discontinuous across the free boundary. Consequently, $u \not \in C^2(\Omega)$ in general. On the other hand, analyzing the free boundary is significantly more complex and remains an active area of research.

\begin{figure}
\centering
  \begin{minipage}[b]{0.48\linewidth}
    \includegraphics[width=\textwidth]{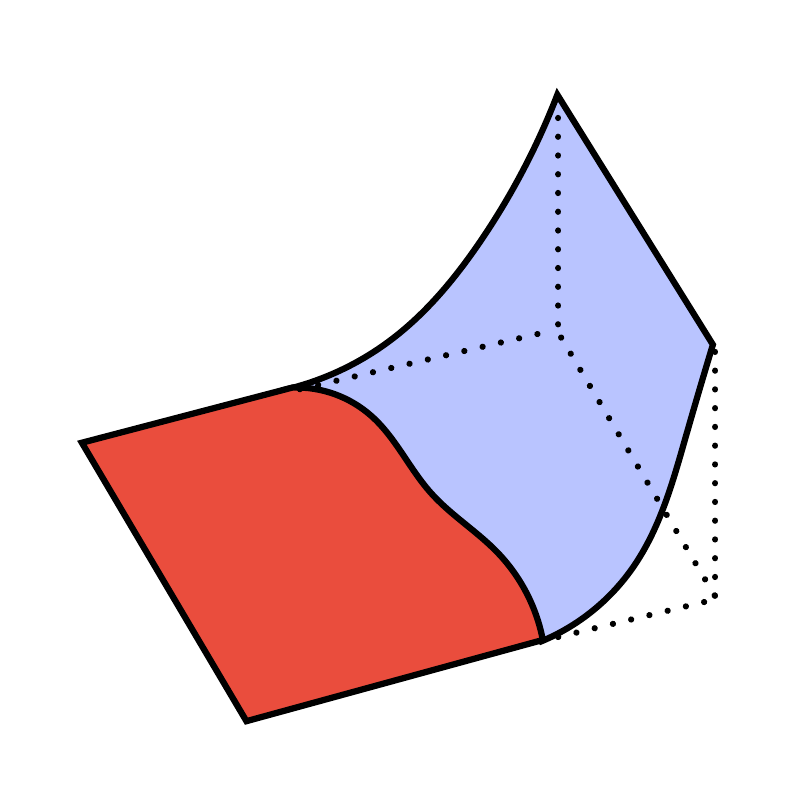}
  \end{minipage}
  \begin{minipage}[b]{0.48\linewidth}
    \includegraphics[width=\textwidth]{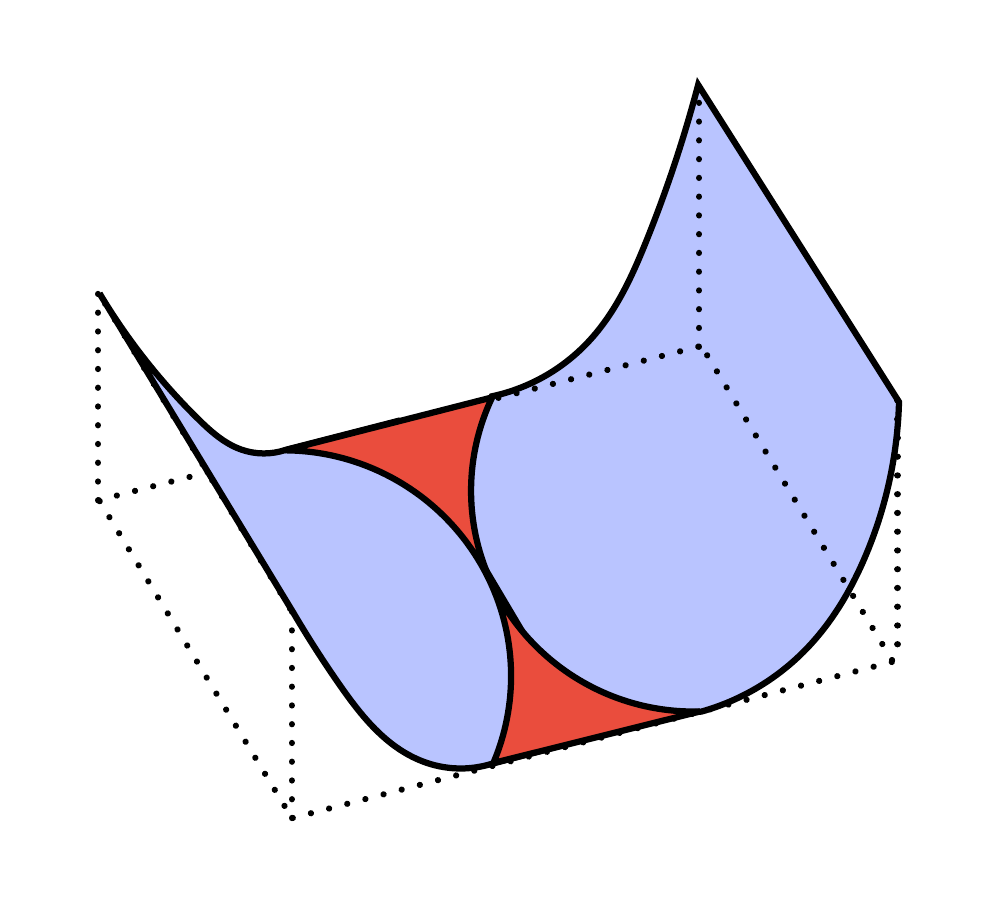}
  \end{minipage}
    \caption{\small Regular vs singular free boundaries: either the contact set looks like a half-space, or it is thin.}
    \label{fig:dichotomy}
\end{figure}

Generally, free boundaries are singular, even when $n = 2$. Caffarelli \cite{Caffarelli1977} showed a powerful dichotomy: when we zoom in around any free boundary point, then
\begin{enumerate}
\item either the contact set looks like a half-space (in which case the case the free boundary is smooth around that point, and we say that the free boundary point is ``regular'');
\item or the contact set is thin (in which case the free boundary point is called ``singular'').
\end{enumerate}
See Figure \ref{fig:dichotomy} for an illustration of either scenario. 

It is well-known that singular points are not rare. Actually, the set of singular points may have the same dimension as the set of regular points.
One can then ask more refined questions about the structure of the set of singular points. 

Caffarelli's result shows that the set of singular points is always \textit{contained} locally in a $C^1$ manifold. Nonetheless,  in specific instances this set can be very wild and display Cantor-like behavior, as was shown by Schaeffer \cite{Schaeffer1977}.   In fact, he conjectured that, although this type of behavior is unavoidable, it is  non-generic: in particular, it should disappear under perturbations of the boundary datum $g$. This conjecture was confirmed first when $n=2$ \cite{Monneau2003} and, more recently, whenever $n\leq 4$ \cite{Figalli2020}. Thus,  in low dimensions, generic free boundaries have no singular points and are therefore smooth. 

\section{The vectorial obstacle problem}

\begin{figure}
\centering
  \begin{minipage}[b]{0.48\linewidth}
    \includegraphics[width=\textwidth]{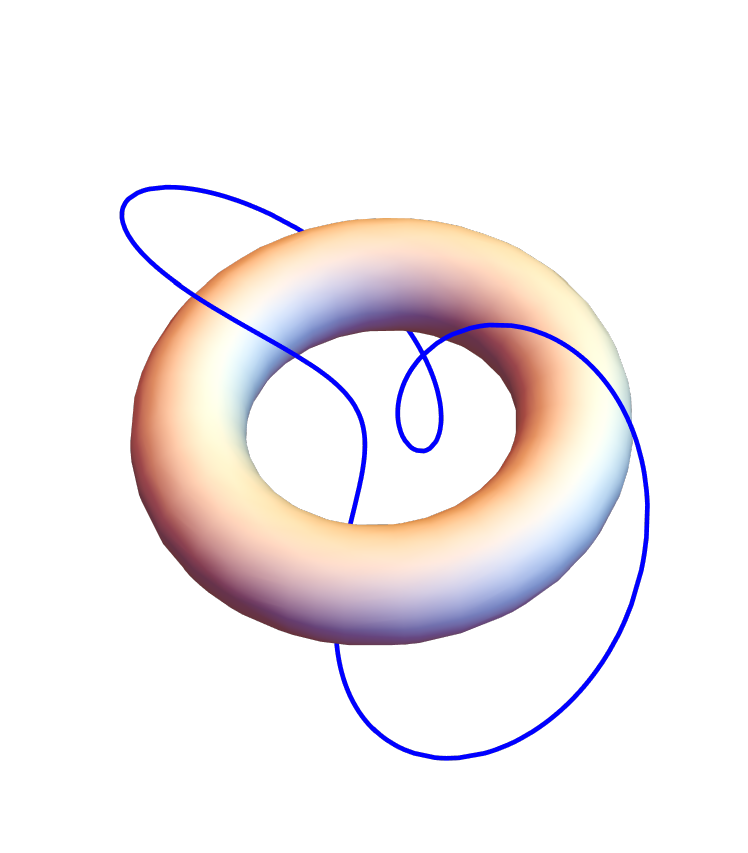}
  \end{minipage}
  \begin{minipage}[b]{0.48\linewidth}
    \includegraphics[width=\textwidth]{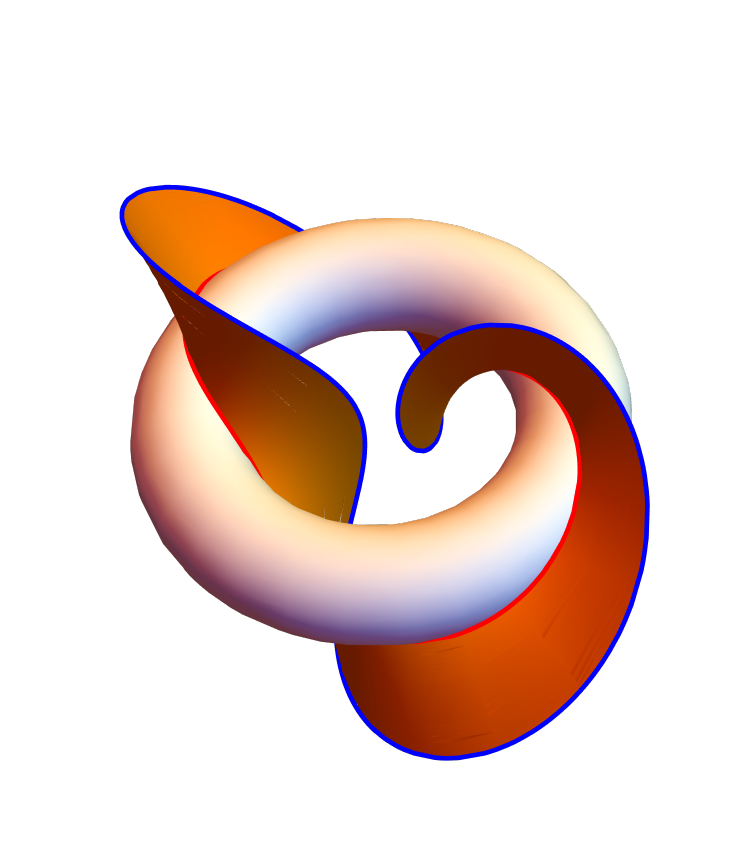}
  \end{minipage}
    \caption{\small Courant's problem: the given surface is a torus and the prescribed Jordan curve is in blue.}
    \label{fig:torus}
\end{figure}

In the previous section we considered graphical deformations of a membrane, see \eqref{eq:graph}. We now consider completely general deformations, which we represent as parametric surfaces with a prescribed boundary. It is natural to consider general deformations since, in more complicated geometric situations, the preferable deformation may not be graphical. For a concrete example, consider the following problem from the classical book of Courant \cite{Courant1977}:
\begin{quote}
{\it ``Find a doubly connected surface of minimal area, one of whose boundaries is free on a given surface, such as a sphere, while the other boundary is a prescribed Jordan curve.''}    
\end{quote}
This is a parametric version of  the King's hat-problem defined in the previous section. Clearly, for complicated Jordan curves, it is unreasonable to expect the solution to be a graph, see e.g.\ Figure \ref{fig:torus}.

In order to formulate precisely the type of problems that we will consider, let $\Omega\subset \R^n, O\subset \R^m$ be smooth bounded domains and let $\g\in C^\infty(\partial \Omega;\R^m)$ be some boundary datum.
As before, we want to minimize
\begin{equation}\label{eq:main}
\cE (\v) = \int_\Omega |D\v|^2\,dx
\end{equation}
among all maps $\v = (v^1, \cdots , v^m)$ that coincide with $\g$ on the boundary and which avoid the obstacle $O$: more precisely, the family of admissible maps is 
$$\widetilde{\K} = \{\v\in W^{1,2}_\g(\Omega;\R^m): \v(x)\not \in O \text{ for a.e.\ } x\in \Omega\},$$
which clearly generalizes our previous definition even in the particular case $m=n+1$. We refer to minimizers of $\cE$ in $\widetilde{\K}$ as \textbf{minimizing constraint maps}.

As in the scalar case, for any admissible map $\v$ the domain $\Omega$ is naturally decomposed into the sets of points which are either mapped to the boundary of the obstacle or to outside the obstacle, together with the free boundary between them:
$$\u^{-1} (\partial O) = \{ x \in \Omega: \u(x)  \in \partial O\},$$
$$\u^{-1} ( \R^m\setminus \overline O) = \{ x \in \Omega:  \u(x) \not \in  \overline O\},$$
$$\cF_\u =  \partial \{ x \in \Omega:  \u(x) \in  \partial O\};$$
see also Figure \ref{Fig:disc} for an illustration.
    Given a minimizer $\u$ of $\mathcal E$  we can then derive equations similar to \eqref{eq:scalarobstacle} to find that 
\begin{equation}\label{eq:main-sys2}
\Delta \u = \cA_\u(D\u,D\u)\,\chi_{\u^{-1}(\partial O)}\quad\text{in }\Omega,
\end{equation}
as was done in \cite{D}.
Here, $\cA$ denotes the second fundamental form of $\partial O$: $\cA_{\u(x)}$  is a quadratic form on the tangent space of $\partial O$ at $\u(x)$ which encodes the curvature of $\partial O$ at that point. For instance, in the model case where $O$ is the  ball $B_1(0)$, we have 
$$\cA_\u(D\u,D\u) = -|D\u|^2 \u.$$ 
Note that \eqref{eq:main-sys2} is a coupled system of $m$ equations, one for each component of $\u$. Moreover,  as in \eqref{eq:scalarobstacle}, we have $\Delta \u=0$ at points where the obstacle is not felt. There is however an important difference compared to \eqref{eq:scalarobstacle}: at points where the obstacle is felt, we have  $\Delta u = \Delta \psi$ in the scalar case, while in the vectorial case we have
\begin{equation}
    \label{eq:harmonicmap}
    \Delta \u = \cA_\u(D\u,D\u).
\end{equation} 
Thus, unlike in the scalar case, the right-hand side in this equation is not fixed a priori but instead depends nonlinearly on the solution itself: in fact, it can be viewed as   a Lagrange multiplier, enforcing the constraint $\u\in \partial O$. Equation \eqref{eq:harmonicmap} asserts that $\u$ is a \textbf{harmonic map} into $\partial O$ when restricted to the interior of the set $\u^{-1}(\partial O)$.

There is an additional  important difference between the scalar and vectorial cases. We saw that, in the scalar case, any minimizer of $\cE$ is a solution of \eqref{eq:scalarobstacle};  conversely, it is known that any solution of \eqref{eq:scalarobstacle} is a minimizer of $\cE$. However, this is not true in the vectorial setting: while minimizers are indeed solutions of \eqref{eq:main-sys2}, the converse  no longer holds. In fact, there are fundamental differences between minima and general solutions of the differential equations.

We now overview some of the history and fundamental examples of minimizing constraint maps:
\begin{enumerate}
    \item When $n=1$, \eqref{eq:main-sys2} becomes a system of ordinary differential equations. In this case, minimizing constraint maps  are simply unit-speed parametrizations of minimizing \textbf{geodesics} in the manifold-with-boundary $\R^m\setminus O$, see Figure \ref{Fig:geodesic}. Note that, even in this case, a solution of the differential equations (i.e., a  geodesic in $\R^m\setminus O$) is not necessarily minimizing.
    \begin{figure}
    \centering
\includegraphics[width=\linewidth]{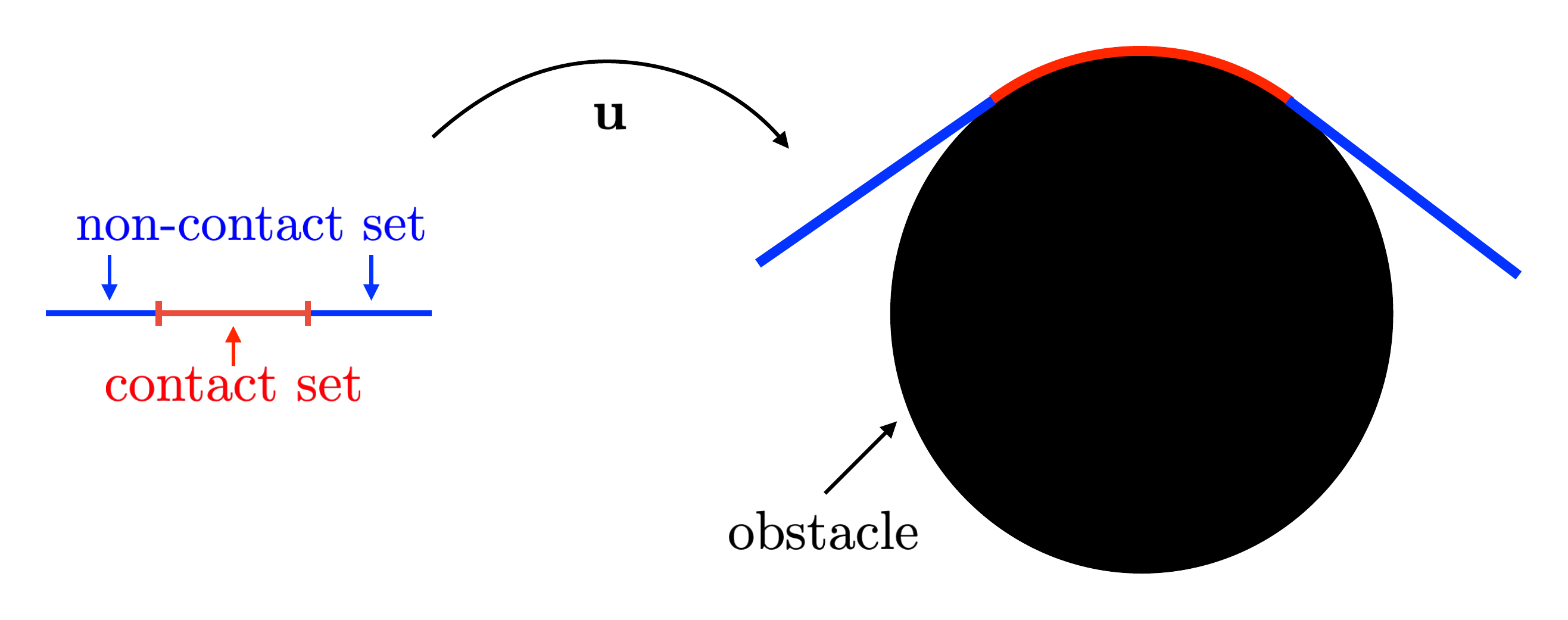}
    \caption {\small A geodesic in $\R^m\setminus O$.} \label{Fig:geodesic}
\end{figure}
    \item The case $n=2$ is very classical and was studied by Morrey \cite{Morrey1948} (for harmonic maps) and then by Hildebrandt \cite{H1} and Tomi \cite{T1} (for constraint maps). Among other things, they proved that minimizing constraint maps are $C^{1,\alpha}$, for any $\alpha\in (0,1)$.
    \item When $n\geq 3$, the theory of minimizing constraint maps becomes especially rich: even in the simplest possible setting where $n=m=3$, $\g=\textbf{id}$ and $\Omega=B_2(0)$, $O=B_1(0)$, minimizing constraint maps develop \textbf{discontinuities}. Indeed, as it is apparent from Figure \ref{Fig:disc}, one has to dissect and break apart   $B_2(0)$ to fit it  onto the target set; this can be made precise using tools from Algebraic Topology and, in particular, degree theory. We shall henceforth  write
    $$\Sigma_\u=\{x\in \Omega: \u \text{ is discontinuous at } x\}.$$
    Starting from the seminal work of Schoen and Uhlenbeck \cite{ScU1}, it was shown in \cite{DF,L,FKS} that 
	\begin{equation}
	\label{eq:dimsing}
	\u\in  C^{1,1}_\text{loc}(\Omega\setminus \Sigma_\u), \qquad \dim\Sigma_\u \leq n-3,
	\end{equation}
    a result which is, in general, optimal. Thus, although discontinuities are inevitable, they form a rather small set. We note that, in special circumstances,  discontinuities may not be present at all: this is the case, for instance,  whenever the obstacle is graphical \cite{FW}. 
    \begin{figure}
    \centering
    \includegraphics[width=\linewidth]{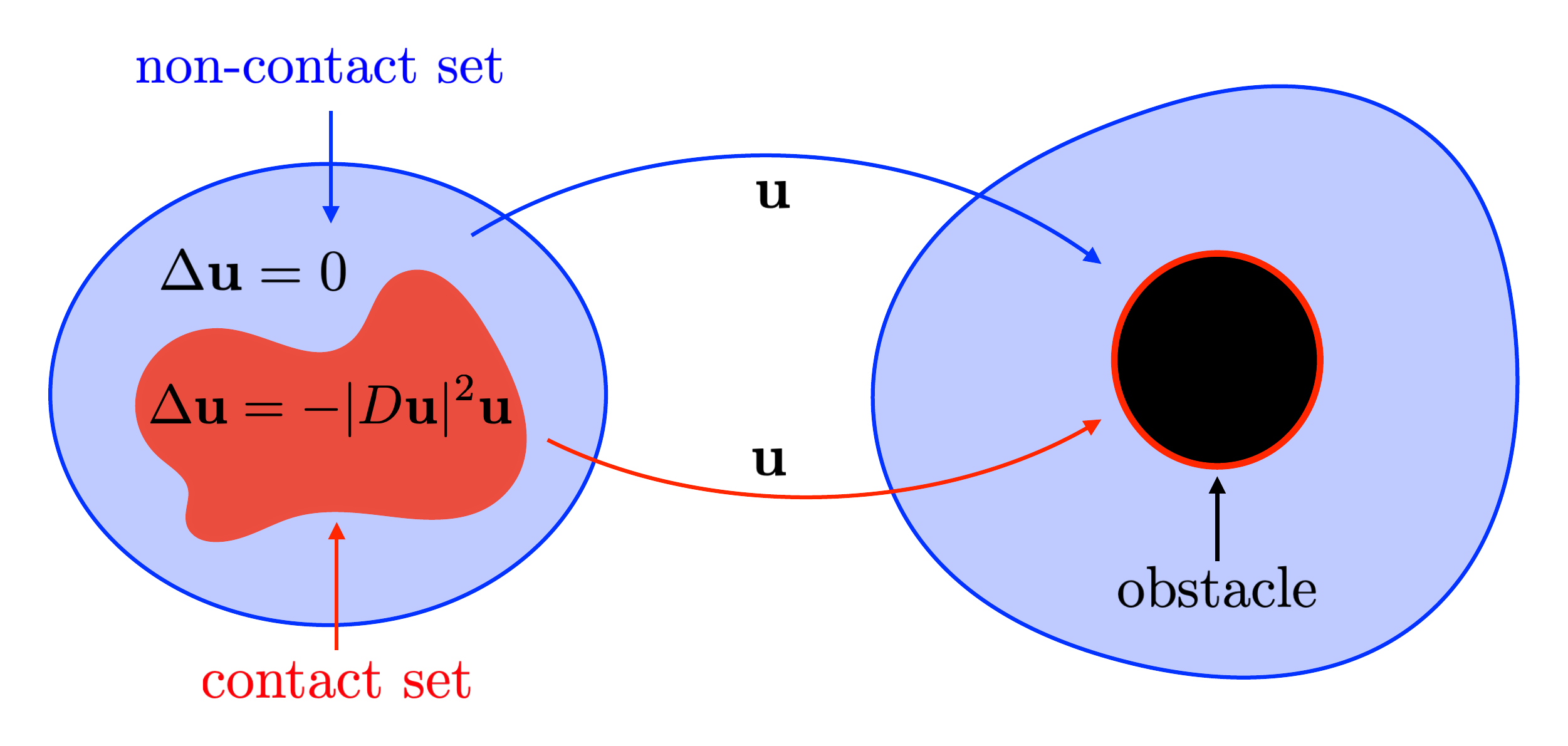}
    \caption {\small An obstacle that forces discontinuities.} \label{Fig:disc}
\end{figure}
\end{enumerate}

It is worth explaining in more detail the tools used in the proof of \eqref{eq:dimsing}.  The main idea, which is standard in this type of problems,  is to understand the behavior of a minimizing constraint map $\u$ around a point $x_0\in \Omega$ by studying rescalings  of $\u$, more precisely
$$\u_{x_0,r}(x):=  \u(x_0+ r x)$$
for $r>0$.  Thus, as $r \searrow 0$, $\u_{x_0,r}$ resembles an increasingly zoomed-in version of $\u$. 
The result \eqref{eq:dimsing}
follows in a standard way from the following three fundamental properties of these rescalings:

\begin{itemize}
\item \textbf{Monotonicity formula:} the rescaled energy
$$E(\u,x_0,r):=r^{2-n} \int_{B_r(x_0)} |D\u|^2 dx$$
is a non-decreasing function of $r$, and it is constant if and only if 
$\u_{x_0,r}(x)= \u(x_0+x)$ for all $x$ and all $r>0$. Note also that $E(\u,x_0,r) = E(\u_{x_0,r},0,1)$, thus the energy is rescaled according to the rescalings of $\u$.
\item \textbf{Compactness:} the sequence of rescalings $\{\u_{x_0,r}(x)\}_{r>0}$ is pre-compact in $W^{1,2}$: if $r_j\searrow 0$ then there is a new subsequence $r_{j'}\searrow 0$ and a minimizing constraint map $\vphi\colon \R^n\to \R^m\setminus O$ such that
$$D \u_{x_0,r_j'} \to D \vphi \quad \text{ in } L^2_{\textup{loc}}(\R^n).$$
By the monotonicity formula, $\vphi$ is 0-homogeneous,  i.e., $\vphi(r x) = \vphi(x)$ for all $x\in \R^n,$ all $r>0$.
\item \textbf{$\e$-regularity theorem:} there is $\e>0$ such that
$$E(\u,x_0,r)\leq \e \quad \implies \quad \u \in C^{1,1}(B_{r/2}(x_0)),$$
with a corresponding estimate.
Thus $\u$ is regular around points of small normalized energy and, conversely,
$$\Sigma_\u=\left\{x_0\in \Omega: \lim_{r\to 0} E(\u,x_0,r)>0\right\}.$$
By the compactness result, we also have $x_0\in \Sigma_\u$ if and only if there is a sequence $r_j\searrow 0$ such that $D \u_{x_0,r_j}\to D\vphi$ in $L^2_\textup{loc}$, for some non-constant map $\vphi$.
\end{itemize}

The results above provide a broad overview of the regularity theory for minimizing constraint maps, paralleling the regularity theory of harmonic maps. However, they leave many fundamental questions unanswered. For instance, how does the geometry of the obstacle $O$ interact with the discontinuity set? And what can be understood about the structure and regularity of the associated free boundaries? We will explore these questions in the following sections.

\section{Discontinuities vs free boundaries}

As we have seen, even in very simple situations, minimizing constraint maps develop discontinuities. Such discontinuities must necessarily lie within the closure of the contact set $\u^{-1}(\partial O)$, since in the non-contact set we have $\Delta \u = 0$ by \eqref{eq:main-sys2}, making $\u$ smooth there. A natural question, then, is whether $\u$ can have discontinuities on the free boundary or if they are confined to the interior of $\u^{-1}(\partial O)$.

The above question is of course only meaningful if a free boundary actually occurs, and this happens only  when $\Omega$ is mapped into a part of $\partial O$ with positive curvature. Indeed, Duzaar \cite{D} showed that if $\nu$ is the outward unit normal to $\partial O$, then 
$$\cA_\u(D\u,D\u) \cdot (\nu\circ \u) \leq 0 \quad \text{ on } \u^{-1}(\partial O).$$
Thus $\u$ may only touch $\partial O$ at points where $\partial O$ looks like the boundary of a convex body.
We leave it to the astute reader to figure out how this inequality results in the existence or non-existence of free boundaries.

The above discussion suggests that we should consider \textbf{convex obstacles} in order to hope for a reasonable answer to our question.  In \cite{FGKS,FGKS2} we showed that indeed the convexity condition leads to some good properties of constraint maps:

\begin{thm}\label{thm:distreg}
    Let $\u\colon \Omega\to \R^m\setminus O$ be a minimizing constraint map, where we assume that 
    $$O \text{ is convex}.$$
    Then the function $x\mapsto \dist(\u(x), O)$ is continuous.
\end{thm}

In particular, Theorem \ref{thm:distreg} guarantees that, for convex obstacles, the non-contact set 
$$\u^{-1}(\R^m\setminus O)=\{\dist(\u, O)>0\}$$ is well-defined and open.

We now provide a sketch of the proof of Theorem \ref{thm:distreg}. 
Since $\dist(\u, O)$ is continuous at points where $\u$ is continuous, it suffices to establish the continuity of $\dist(\u, O)$ at the discontinuity points of $\u$. For any such point, say $x_0 \in \Sigma_\u$, as shown in the previous section, there exists a non-constant, minimizing constraint map $\vphi$ that is 0-homogeneous and arises as the limit of a sequence of rescalings $\u_{x_0, r}$.
As $\vphi$ is itself a solution of the Euler--Lagrange system \eqref{eq:main-sys2} and, as $O$ is convex,  we have that the function
$$x\mapsto \dist(\vphi(x),O) \quad \textup{ is subharmonic}.$$
Since this function is also 0-homogeneous, by the maximum principle it must be constant. In fact, we must have
$$\dist(\vphi,O)=0,$$
for otherwise $\vphi$ would never touch $\partial O$ and so we would have $\Delta \vphi=0$, which is impossible since $\vphi$ is 0-homogeneous and non-constant. Now, according to the definition of the rescalings, we have
\begin{align*}
\|\dist(\u,O)\|_{L^2(B_r(x_0))} & = \|\dist(\u_{x_0,r},O)\|_{L^2(B_1(0))} \\
&\to \|\dist(\vphi,O)\|_{L^2(B_1(0))} = 0,
\end{align*}
as $r\to0$,  since $\u_{x_0,r}\to \vphi$ in $L^2_\textup{loc}$. To conclude that $\dist(\u, O)$ is continuous at $x_0$, it suffices to upgrade the above $L^2$-convergence to $L^\infty$-convergence. Since $\dist(\u, O)$ is subharmonic (due to the convexity of $O$), the $L^2$-norm controls the $L^\infty$-norm by the mean value property.

\medskip

Note that we may always decompose a map $\u\colon \Omega\to\R^m\setminus O$ in its normal and tangential components (with respect to $O$); explicitly, we have
$$\u = \dist(\u,O)(\nu \circ \Pi \circ \u) + \Pi\circ \u,$$
where as before $\nu$ is the outwards unit normal to $O$ and $\Pi\colon \R^m\setminus O\to \partial O$ is the nearest point projection.
Theorem \ref{thm:distreg} shows that the components of $\u$ normal to $O$ are continuous, but what about the tangential components? Surprisingly, in  \cite{FGKS2} we observed that the tangential components are in general discontinuous on the free boundary:

\begin{thm}\label{thm:flat}
Let $n = m$ and $\Omega = B$ be a ball. There exists a convex obstacle $O$ such that the energy-minimizing map $\u \colon \Omega \to \R^n \setminus O$ has the property that
$$\Sigma_\u \cap \cF_\mathbf u \neq \emptyset.$$

\end{thm}

\begin{figure}
    \centering
    \includegraphics[width=\linewidth]{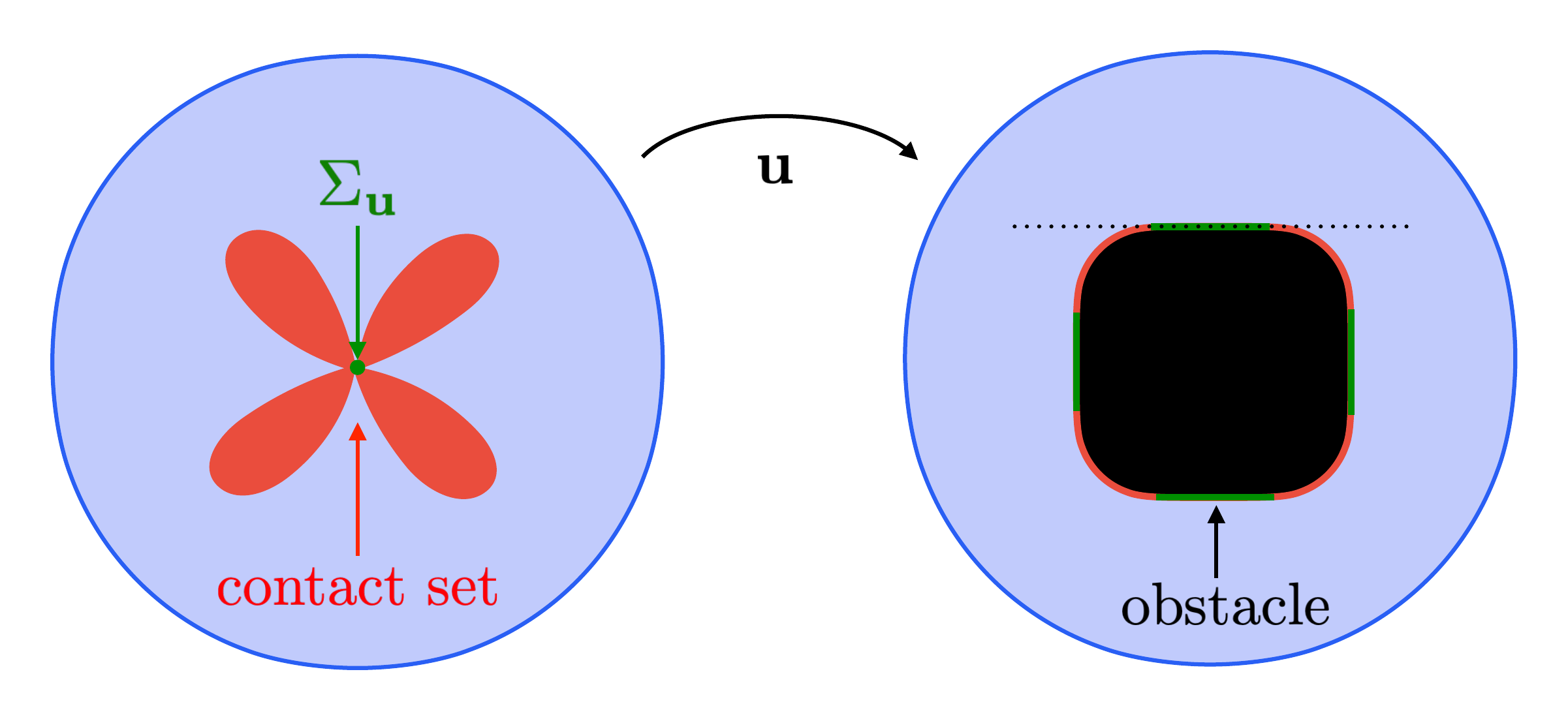}
    \caption {\small A convex obstacle with flat sides forces discontinuities on the free boundary: the discontinuity points are mapped to the flat sides.} \label{Fig:flat}
\end{figure}

The obstacles in Theorem \ref{thm:flat} are easy to describe:  any convex obstacle with \textbf{flat pieces} (i.e. an obstacle which contains a portion of a hyperplane in $\R^m$) will induce discontinuities on the corresponding free boundaries, provided that 
\begin{equation}
\label{eq:idbc}
\overline O \subset  \Omega=B, \qquad \u = \textbf{id} \text{ on } \partial \Omega,
\end{equation}
see Figure \ref{Fig:flat}. This is somewhat counter-intuitive: when the obstacle is a half-space we have $\cA=0$; so minimizing constraint maps are simply solutions of $\Delta \u=0$ and, in particular, they are analytic. 
 
We now explain the mechanism behind Theorem \ref{thm:flat}.  Suppose\footnote{Of course, the main challenge here is to prove the existence of a free boundary point that is mapped to a flat piece. We overcome this issue in \cite{FGKS2} by a topological degree argument combined with some tools in PDEs and differential geometry.} there is a free boundary point $x_0\in \cF_\u$ that is a continuity point of $\u$ and which is mapped to a flat piece $P\subset \partial O$.  Any such flat piece can be  characterized by a linear equation; for simplicity, let us say  
$$P=\{\mathbf y=(y^1,\ldots,y^n): y^n=\tfrac 1 2 \}\cap O, $$ 
 and  assume that $\{y^n\geq \tfrac 1 2\} \subset \R^n\setminus O$, 
see Figure \ref{Fig:flat}.  By continuity of $\u$ at $x_0$, we can find a small radius $r>0$ such that $$\u(B_r(x_0) \cap \u^{-1}(\partial O))\subset P.$$ Since $\cA = 0$ in $P$, it follows from  \eqref{eq:main-sys2} that  $\Delta \u = 0$ in $B_r(x_0)$. Thus 
$$
\Delta u^n=0 \text{ and } 
u^n\geq \tfrac 1 2 \text{ in } B_r(x_0),\quad u^n(x_0)=\tfrac 1 2.
$$
The maximum principle for harmonic functions forces $u^n=\tfrac 1 2$ in $B_r(x_0)$. Since $\Delta \u =0$ in the non-contact set $\u^{-1}(\R^n\setminus O)$ and the latter is connected, the unique continuation property of harmonic functions implies that $u^n = \tfrac{1}{2}$ throughout the entire non-contact set (recall that harmonic functions are analytic). This, of course, contradicts our choice of boundary conditions in \eqref{eq:idbc}.
Thus, $\u$ cannot  continuously map  a free boundary point into a flat piece.  On the other hand,  \eqref{eq:idbc} requires  that  $\u$ to fill the annulus $\Omega\setminus O$, i.e.
$$\u(\Omega)=\Omega\setminus \overline O,$$ in order to minimize the amount of  ``breaking'' occurring in $\Omega$. Consequently, there must be free boundary points that map to a flat piece, and these points are necessarily discontinuities.

It is an interesting problem to understand the geometry of the contact set $\u^{-1}(\partial O)$ for the above example. We refer the reader to Figure \ref{Fig:flat} for a suggestion of what it may look like; we emphasize however that we do not yet  know  if the free boundary is singular at points in $\Sigma_\u \cap \cF_\u$. 
The  argument above  suggests that  $\u$ should map a single free boundary point to an entire  flat piece of $\partial O$.  Therefore,  there may be a singular point in $\cF_\u$ for each flat piece of $\partial O$, or  there could be a singular point in $\cF_\u$ that  is mapped simultaneously to different flat pieces.

\medskip
Theorems \ref{thm:distreg} and \ref{thm:flat} give a rather complete picture of what happens for general convex obstacles. 
However, in the scalar obstacle problem, where $\cA$ is replaced by $\Delta \psi$ for some prescribed function $\psi$, one typically assumes the quantitative condition \eqref{eq:Deltapsi=1}. Thus it is natural to wonder whether a more quantitative version of convexity, such as \textbf{uniform convexity}, leads to better results.  We recall that $O$ is said to be uniformly convex if all of the principal curvatures of $\partial O$ are positive and bounded away from zero: essentially, $O$ looks ``round''. In \cite{FGKS2}, we obtained the following positive result: 

\begin{thm}\label{thm:unireg}
Let $\u\colon \Omega\to \R^m\setminus O$ be a minimizing constraint map, where we assume that 
$$O \text{ is \textbf{uniformly} convex.}$$
There exists  $\delta>0$ such that
$$x\in \Sigma_\u \quad \implies \quad \delta \leq \dist(x,\cF_\u).$$
In other words, that are no discontinuities on the free boundary.
\end{thm}

The proof of Theorem \ref{thm:unireg} is unfortunately too technically involved  to describe in detail here.
 In essence, the uniform convexity of $O$ allows for a quantification of the proof of Theorem \ref{thm:distreg} described above, but it requires many new ideas.  The proof integrates concepts from De Giorgi regarind  subsolutions of elliptic equations,  the frequency function introduced by Almgren in the study of minimal surfaces, and the identification of a critical scale which separates precisely the regions of 	regular and irregular behavior of $\u$, in line with   the recent work of Naber and Valtorta \cite{Naber2017}.

\section{Branch points and singular free boundaries}

In the scalar obstacle problem,  a successful  regularity theory for free boundaries was built under the quantitative assumption \eqref{eq:Deltapsi=1}. The scalar theory can be applied to study the regularity of free boundaries also in the vectorial setting, provided that the right-hand side in \eqref{eq:main-sys2} is non-degenerate, the relevant quantity here being
$$F(\u)(x):=\Hess \left(\dist(\cdot,O)\right)_{\u(x)}[D\u,D\u]\geq 0,$$
where the last inequality holds whenever $O$ is convex.
It is not surprising that the Hessian of the distance to $O$ is relevant here, since such a quantity is closely related to the second fundamental form of $\partial O$.

Essentially, whenever $F(\u)>0$, condition \eqref{eq:Deltapsi=1} holds and the scalar theory is applicable \cite{FKS}.  There are two important cases, however, where this condition cannot hold:
\begin{enumerate}
\item at discontinuity points of $\u$, then $|F(\u)|=\infty$;
\item at points where $D\u=0$, then $F(\u)=0$.
\end{enumerate}
Thus, in both cases, the scalar theory is entirely inapplicable. 
In the previous section, we provided a comprehensive overview  of how the geometry of the obstacle interacts  with the existence of discontinuity points on free boundaries. We will now turn our attention to discussing the set of points  where 
$$D\u=0.$$
We refer to that set as the set of \textbf{branch points}. 

There are two  distinct yet  intriguing aspects of branch points. Firstly, if branch points exist on a free boundary, they may significantly affect its regularity, as the vectorial obstacle problem becomes degenerate at these points. We will discuss this aspect in detail shortly. Secondly, at a branch point, one cannot apply the Implicit Function Theorem to $\u$. Thus, if we consider $\u$ as a parametrization of its image $\u(\Omega)$, the image itself may develop singularities.

\begin{figure}[h]
\centering
  \begin{minipage}[b]{0.48\linewidth}
    \includegraphics[width=\textwidth]{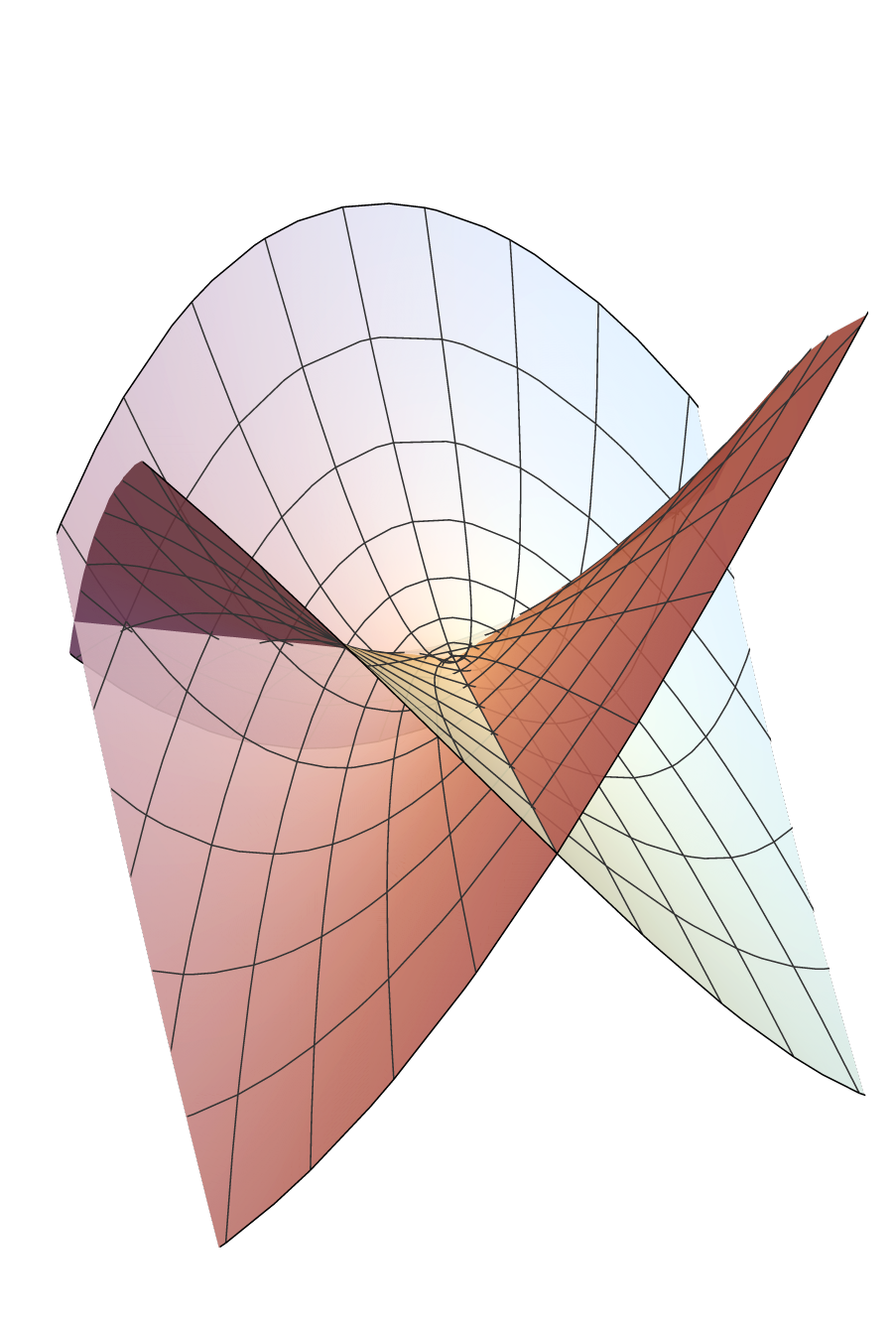}
  \end{minipage}
  \begin{minipage}[b]{0.48\linewidth}
    \includegraphics[width=\textwidth]{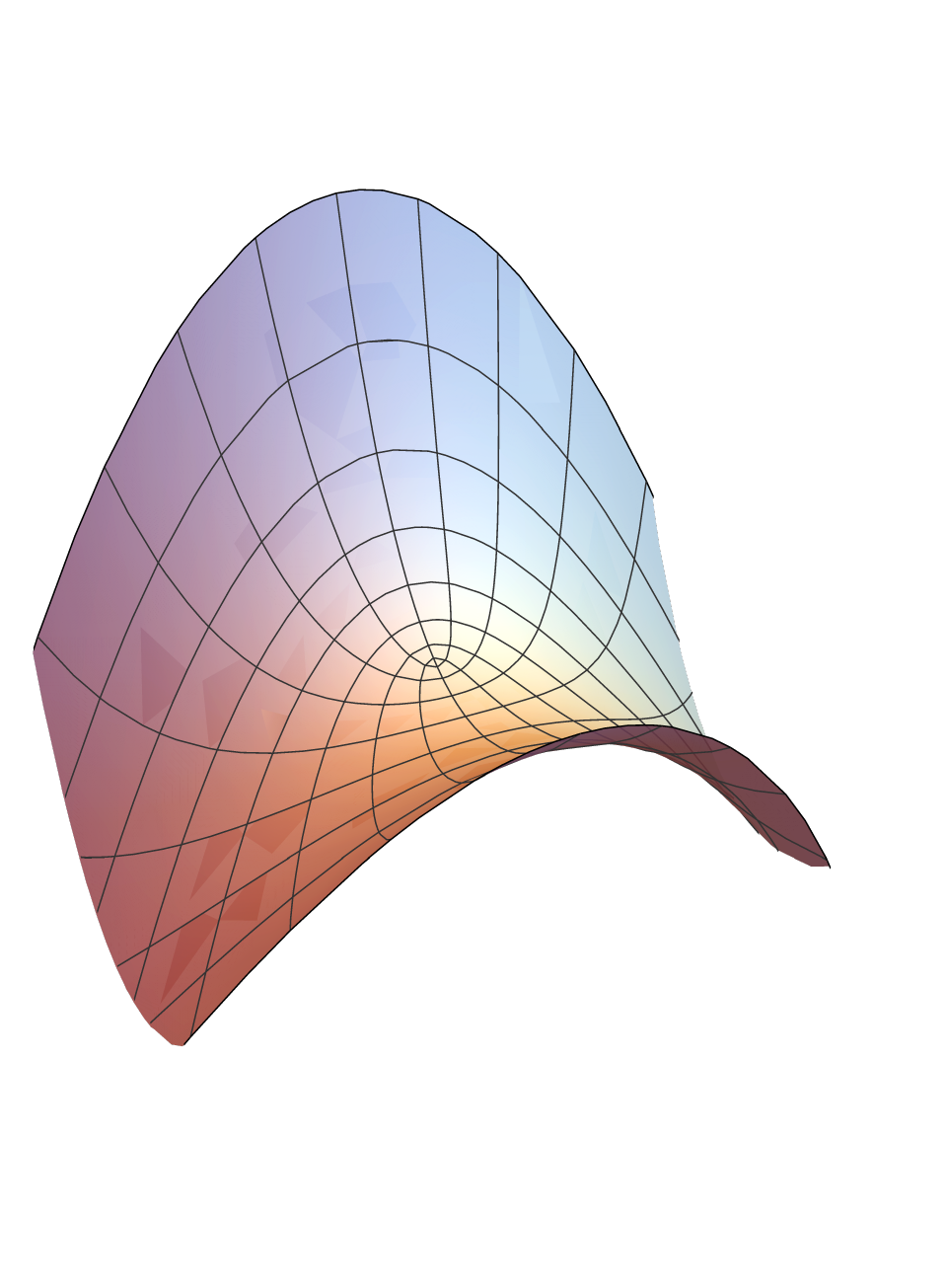}
  \end{minipage}
    \caption{\small $\u_k(\R^2)$ for $k=3$ (left) and $k=4$ (right).}
    \label{fig:branch}
\end{figure}

Let us illustrate this second point in a simple but compelling example, which will also clarify how our notion of branch point is related to analogue concepts in classical function theory and in the theory of parametric minimal surfaces.  Consider, for $k \in \N$, the maps $\u_k\colon \C\cong \R^2\to \R^3$ defined by 
$$\u_k(z) = (z^2, \textup{Re}\, z^k),  \quad z\in \C,$$
whose images are depicted in Figure \ref{fig:branch}.  
Notice that 0 is a branch point of $\u_k$ (and it is, in fact, the only such point). There is an interesting difference between the case where $k$ is odd and $k$ is even: in the former, there is an actual singularity of $\u(\R^2)$ at 0; in the latter, the parametrization $\u$ is degenerate at 0, but the surface $\u(\R^2)$ is perfectly smooth. 

\begin{figure}
    \centering
    \includegraphics[width=\linewidth]{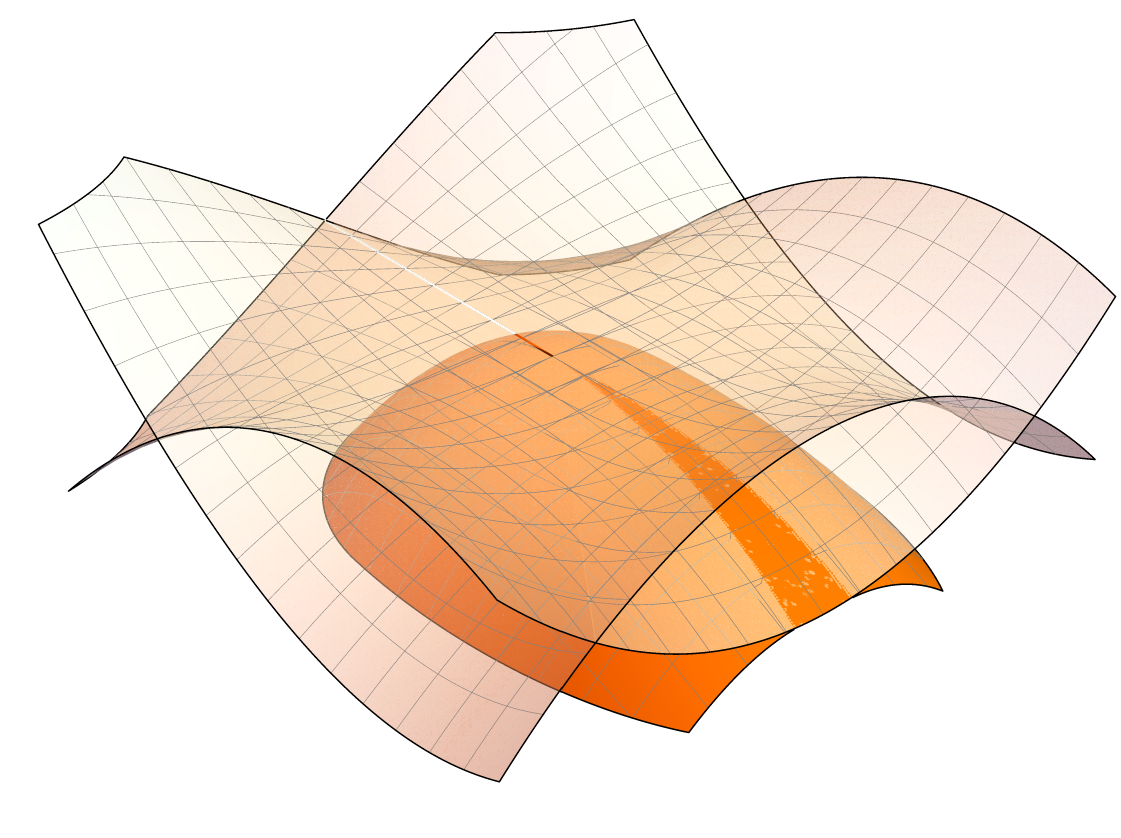}
    \caption {\small An artificial obstacle touching the image of $\u_5$ over a cone with tip at 0.  Note that the obstacle is not convex but that it is of class $C^{2,1/2}$.
} \label{fig:artificial}
\end{figure}

The maps $\u_k$ clearly satisfy $\Delta \u_k=0$, and so they are minimizers of the Dirichlet energy \eqref{eq:main},  without any constraints on the image of the competitors. In particular, they are also minimizing constraint maps for any obstacle $O\subset \R^3$ such that 
$$O\subset \R^3 \setminus \overline{\u_k(\R^2)}.$$  By choosing appropriate sets $O$, one can produce simple but interesting examples of branch points on free boundaries,  as the reader can see in Figure \ref{fig:artificial}.  The case where $k$ is odd is especially interesting: as in Figure \ref{fig:artificial},  for any such $k$ one can construct a non-convex obstacle of class $C^{(k-1)/2}$ such that the free boundary is a cone with tip at the branch point at 0. 

The type of free boundary singularities depicted in Figure \ref{fig:artificial} cannot occur in the scalar obstacle problem under the quantitative condition \eqref{eq:Deltapsi=1},  as per Caffarelli's dichotomy. As the above obstacle is not even convex, and in light of the positive results of Theorem \ref{thm:unireg}, it is natural to ask what happens for uniformly convex obstacles: can the free boundary prevent branching?

It turns out that, even in very simple situations, the answer is no, even when the obstacle is the unit ball $\mb B^3\subset \R^3$:

\begin{thm}\label{thm:branch}
There is a minimizing constraint map $\u\colon \Omega\to \R^3\setminus \mb B^3$ such that 
$$\{D\u=0\}\cap \cF_\u \neq \emptyset,$$ and the free boundary $\cF_\u$ is singular at such points.
\end{thm}

The geometric idea behind Theorem \ref{thm:branch} is simple to explain. In doing so, we will also provide a more detailed version of the result, including a more precise description of the nature of the singularities of the free boundary.

\begin{figure}
    \centering
    \includegraphics[width=\linewidth]{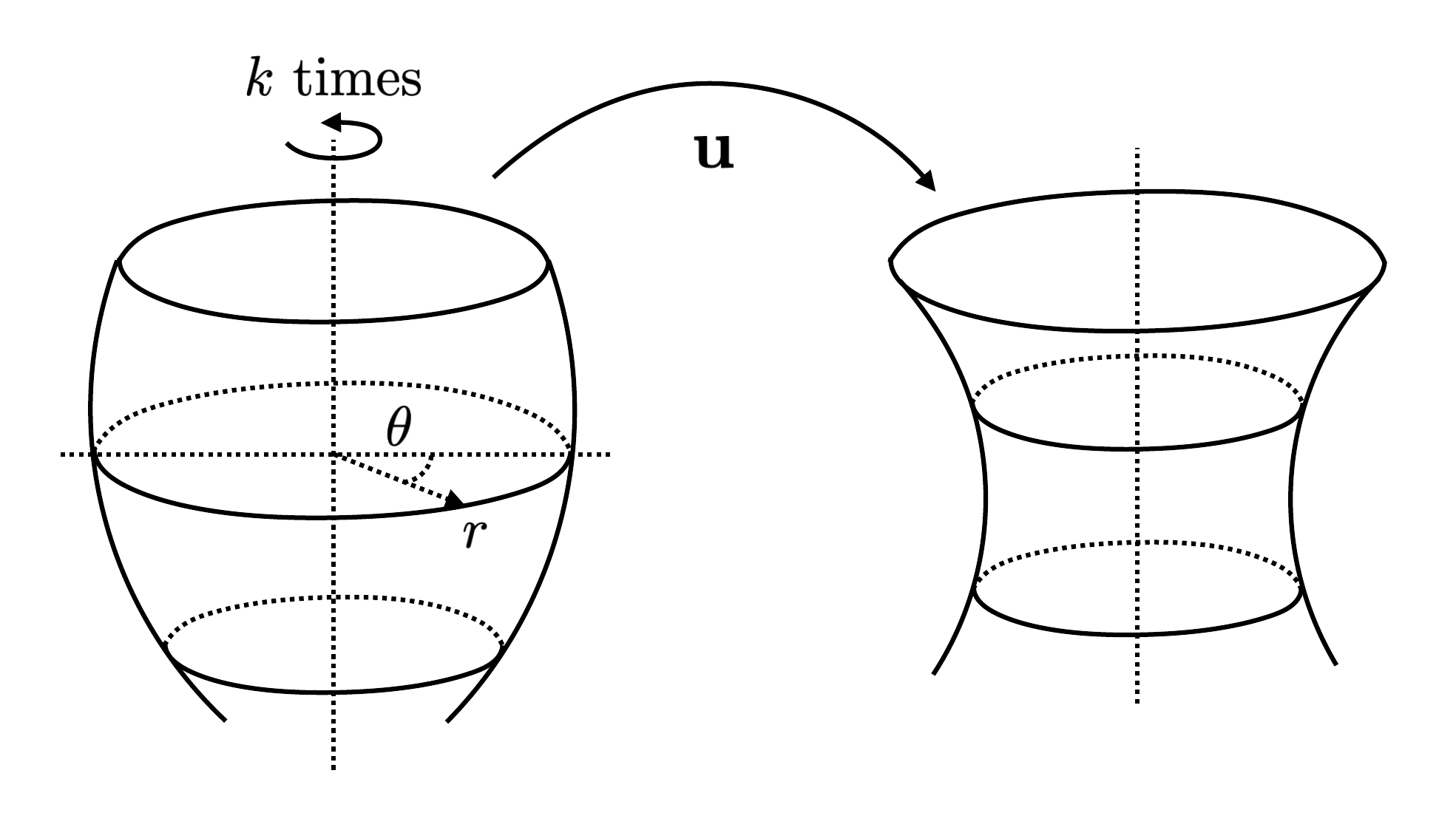}
    \caption {\small A $k$-axially symmetric map $\u$ maps horizontal circles centered on the vertical axis to other such circles, while rotating each of them $k$-times.
} \label{fig:axial}
\end{figure}

For a fixed integer $k\in \N$, we consider $k$-axially symmetric maps $\u\colon \mb B^3\to \R^3\setminus \mb B^3$. Specifically, in cylindrical coordinates $(r,\theta,x_3)$ and $\u=(u^r, u^\theta, u^3)$ for both the domain and the target, we have
$$
\u(r,\theta,x_3) = \left(u^r(r,x_3),  k \theta, u^3(r,x_3)\right).
$$
In other words, $\u$ maps each horizontal circle centered on the vertical axis to another such circle, while rotating it $k$ times, see Figure \ref{fig:axial}. Given any non-trivial $k$-axially symmetric boundary condition, for instance 
$$\g(r,\theta,x_3) = (r,k \theta,x_3),$$
we then minimize the Dirichlet energy \eqref{eq:main} among $k$-axially symmetric maps taking values in $\R^3\setminus \mb B^3$  which agree with $g$ on $\partial \mb B^3$.  Although the resulting map $\u$ is not necessarily a minimizing constraint map, it can be shown that it is in fact minimizing around any of its continuity points. In particular, by (a variant of) Theorem \ref{thm:unireg},  $\u$ is minimizing in a neighborhood of any free boundary point, which is the region that concerns us. 

Let us now take a point on the vertical axis. Our construction is based on the following simple observation: if $k\geq 2$ and $\u=(u^1,u^2,u^3)$ is a $k$-axially symmetric map,  then
\begin{equation}
\label{eq:horizontalcomps}
Du^1 (0,0,x_3)= Du^2 (0,0,x_3)=0,
\end{equation}
provided that $\u$ is $C^1$ around $(0,0,x_3)$.  Indeed, since $\u$ rotates  at least twice about the vertical axis, the only way it can do so in a regular manner  is by having vanishing gradient on the this axis.  Now suppose, in addition, that $(0,0,x_3)\in \cF_\u$ is a free boundary point. Then $|\u|(0,0,x_3)=1$, and  since $|\u|\geq 1$ and $|\u|$ is smooth around $(0,0,x_3)$, we have
$$0 = D|\u| (0,0,x_3) =  D u^3(0,0,x_3),$$
where we used also \eqref{eq:horizontalcomps} in the last equality. We conclude that \textit{any free boundary point on the vertical axis is a branch point}.

\begin{figure}
    \centering
    \includegraphics[width=\linewidth]{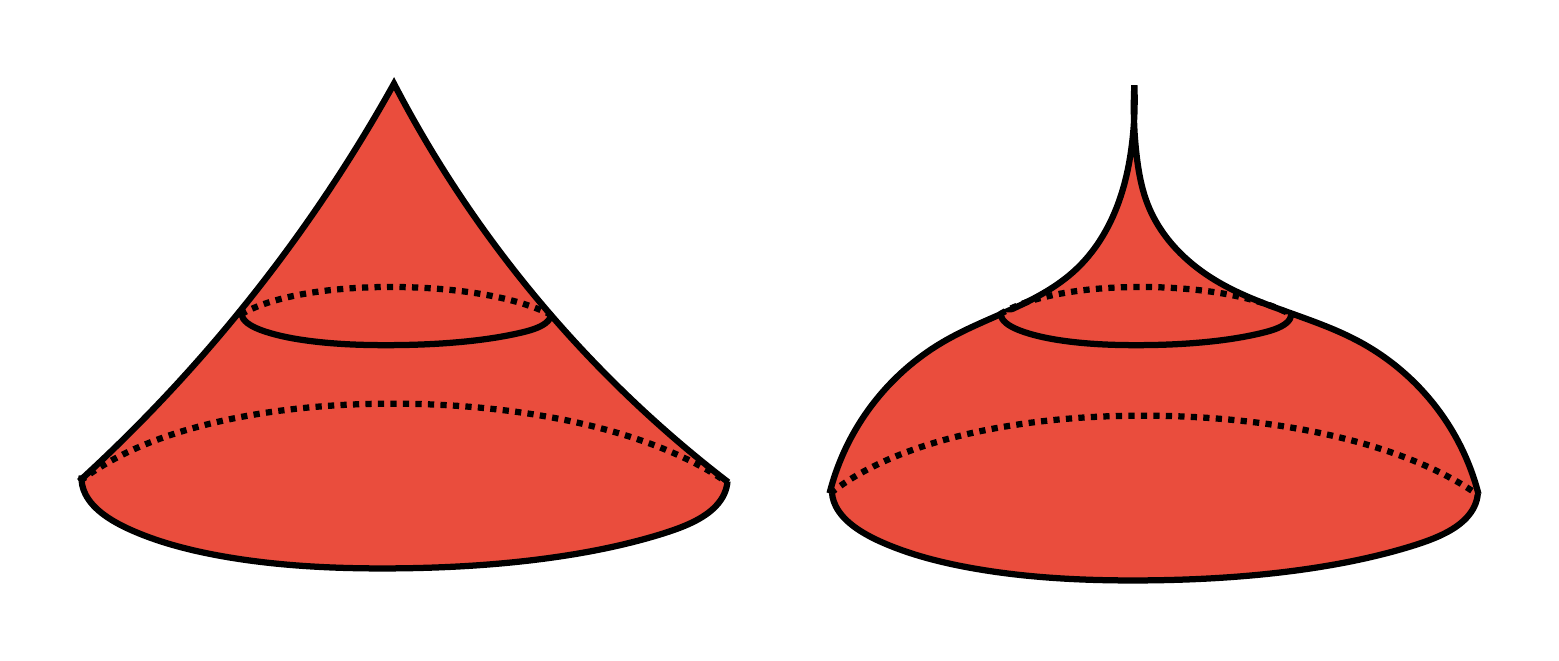}
    \caption {\small At a branch point, the contact set either looks like a cone or like the vertical axis. In either case, the free boundary is singular.
} \label{fig:bp}
\end{figure}

What does the free boundary look like near such a point? Since we are in the axially symmetric setting, we can compute many things fairly explicitly.  In \cite{FGKS2}, we showed that, as we zoom in around a free boundary point on the vertical axis, the free boundary either looks like a cone or like the vertical axis, see Figure \ref{fig:bp}. In the former case, the angle of the cone is quantized: if $\varphi$ is the angle between the cone and the vertical axis, then $\cos(\varphi)$ must be a zero of the Legendre polynomial of order $2k-1$, where $k$ is the number of times $\u$ wraps around the vertical axis. In either case, we see that the free boundary is singular at such a point. 	An interesting problem is to understand whether both types of singularities are actually possible: in our analysis we can only show that the free boundary looks like one of the two scenarios in Figure \ref{fig:bp}, but we do not know whether both may occur in practice. 
Recall that cone-like singularities cannot occur in the scalar obstacle problem (cf.\ Figure \ref{fig:dichotomy}). If they occur in this context, they represent a genuinely vectorial feature resulting from the emergence  of branch points.

\subsection*{Acknowledgements}

A.F. has been partially supported by the European Research Council under the Grant Agreement No.~721675 ``Regularity and Stability in Partial Differential Equations (RSPDE)''.
A.G. was supported by Dr.~Max Rössler, the Walter Haefner Foundation, and the ETH Zürich Foundation. 
S.K. was supported by a grant from the Verg Foundation.
H.S. was supported by the Swedish Research Council (grant no.~2021-03700).

\section*{Declarations}

\noindent {\bf  Data availability statement:} All data needed are contained in the manuscript.

\medskip

\noindent {\bf  Funding and/or Conflicts of interests/Competing interests:} The authors declare that there are no financial, competing or conflict of interests.

\bibliography{references.bib}


\end{document}